\documentclass[12pt]{article}
\usepackage{amsfonts}
\usepackage{amssymb}
\usepackage{enumerate}
\newtheorem{st}      {Theorem}
\newtheorem{prop}  {Proposition}
\newtheorem{lem} {Lemma}
\newtheorem{cor} {Corollary}

\newtheorem{exam}{Example}[section]
\newtheorem{num}  {} [section]

\newcommand{\rmap}{\longrightarrow}
\newcommand{\Boxe}{\raisebox{.8ex}{\framebox}}
\newcommand{\lmap}{\longleftarrow}

\newcommand{\U}{\ensuremath{\mathcal{U}}}
\newcommand{\oo}{\ensuremath{\mathcal{O}}}

\newcommand{\A}{\ensuremath{\mathcal{A}}}

\newcommand{\el}{\ensuremath{\mathcal{L}}}
\newcommand{\F}{\ensuremath{\mathcal{F}}}
\newcommand{\es}{\ensuremath{\mathcal{S}}}

\newcommand{\G}{\ensuremath{\mathcal{G}}}

\newcommand{\ps}{{\raise 1pt\hbox{\tiny (}}}

\newcommand{\pss}{{\raise 1pt\hbox{\tiny [}}}
\newcommand{\pdd}{{\raise 1pt\hbox{\tiny ]}}}
\newcommand{\pd}{{\raise 1pt\hbox{\tiny )}}}

\newcommand{\bs}{{\raise 1pt\hbox{\tiny [}}}
\newcommand{\bd}{{\raise 1pt\hbox{\tiny ]}}}
\newcommand{\nG}[1]{\ensuremath{\G^{\ps #1\pd}}}

\def\cross{\mathinner{\mathrel{\raise0.8pt\hbox{$\scriptstyle>$}}
                 \joinrel\mathrel\triangleleft}}

\def\compose{{\raise 1pt\hbox{$\scriptscriptstyle\circ$}}}
\def\dcross{{\raise 0.5pt\hbox{$\scriptscriptstyle\boxtime$}}}

\advance\textheight by \topskip
\usepackage[all]{xy}
\usepackage{epsf}
\usepackage{amsfonts}
\usepackage{amssymb}

\begin{document}

\title{$\check{C}$ech-De Rham theory for leaf spaces of foliations\thanks{Research supported by NWO}}
\author {Marius Crainic and Ieke Moerdijk}
\pagestyle{myheadings}
\date{}
\maketitle
%\begin{abstract}

%{\it Keywords}:  foliations, characteristic classes, $\check{C}$ech-DeRham complex, \'etale groupoids
%\end{abstract}
%********

%\tableofcontents

%\newpage

%***************************************************************
%***************************************************************
%***************************************************************
%***************************************************************
\section*{Introduction}
%***************************************************************
%***************************************************************
%***************************************************************
%***************************************************************

This paper is concerned with characteristic classes in the cohomology of leaf spaces of
foliations. For a manifold $M$ equipped with a foliation $\F$ it is well-known that the coarse
(naive) leaf space $M/\F$, obtained from $M$ by identifying each leaf to a point, contains very little information.
In the literature, various models for a finer leaf space $M/\F$
are used for defining its cohomology. For example, one considers the cohomology of the classifying space of the foliation
\cite{Bott, Dupon, Haefl, LaPa},
the sheaf cohomology of its holonomy groupoid \cite{CrMo, difcoh, conj}, or the cyclic cohomology of its convolution algebra \cite{Co3, Cra}. 
Each of these methods has considerable drawbacks. E.g. they all involve non-Hausdorff spaces in an essential way.
More specifically, the classifying space, which is probably the most common model for the ``fine'' leaf space,
is a space which in general is infinite dimensional and non-Hausdorff, it is not a CW-complex, and it has lost all the
smooth structure of the original foliation. In particular, it is not suitable for constructing cohomology theories
with compact support. 
For this reason, the construction of characteristic classes in the cohomology of the classifying space of the foliation proceeds
in a very indirect way, and many of the standard geometrical constructions have to  be replaced by or supplied with abstract non-trivial
arguments. The same applies to the construction of ``universal'' characteristic classes in the cohomology of the classifying space of
the Haefliger groupoid $\Gamma^q$. It is possible to construct interesting classes of (foliated or transversal) bundles
over foliations by explicit geometrical methods \cite{Bott, KT}, but these classes are constructed in the cohomology
of the manifold $M$ rather than that of the leaf space $M/\F$.\\
\hspace*{.3in}The purpose of this paper is to present a ``$\check{C}$ech-De Rham'' model for the cohomology of leaf spaces (Section \ref{CDRcomplex}), which circumvents
the problems mentioned above. This $\check{C}$ech-De Rham model lends itself to the construction of (known) characteristic classes,
now by explicit geometrical constructions which are immediate extensions of the standard constructions for manifolds (Section \ref{classes}).
As a consequence, for any transversal principal bundle over a foliated manifold $(M, \F)$, we are able to lift the characteristic classes
constructed in $H^*(M)$ by the methods of \cite{KaTo}, to the $\check{C}$ech-De Rham cohomology $H^*(M/\F)$, and establish all the relations, 
such as the Bott vanishing theorem, at the level of $H^*(M/\F)$ (see Theorem \ref{theorem2} below).\\
\hspace*{.3in}We want to emphasize that the construction of the $\check{C}$ech-De Rham model and of the characteristic classes makes no reference
to (holonomy) groupoids or classifying spaces. In particular, there are no non-Hausdorffness problems, and these constructions can 
be understood by anyone having some background in differential geometry, including familiarity with the very basic definitions
concerning foliations.\\ 
\hspace*{.3in}To prove that our $\check{C}$ech-De Rham model gives in fact the same cohomology as the other models (Theorem \ref{theorem1}), we use \'etale groupoids (Section \ref{etale}). 
In fact, our model and the associated method for constructing characteristic classes applies to any \'etale groupoid, not just to holonomy groupoids
(see Theorem \ref{theoreml}, and \ref{CWetalgr}).
In particular, when used in the context of the Haefliger groupoid $\Gamma^q$, it provides an explicit geometric construction of the
universal geometrical characteristic classes (as a map from Gelfand-Fuchs cohomology into the cohomology of $B\Gamma^q$ \cite{Bott}). 
In this way we rediscover (and explain) the Thurston formula and the Bott formulas for cocycles on diffeomorphism groups \cite{Boform}
(for these explicit formulas, see Section \ref{explicit}).
Other groupoids of interest, different from holonomy groupoids, are the monodromy groupoids of foliations. Our methods also
show that the characteristic classes of foliated bundles \cite{KaTo} actually live in the cohomology of the monodromy
groupoid of the foliation, rather in the cohomology of $M$ itself. \\
\hspace*{.3in}Our $\check{C}$ech-De Rham cohomology also has a natural version with compact supports, which is related to the one 
with arbitrary supports by an obvious duality. When passing to the cohomology of holonomy groupoids, this duality becomes the 
Poincar\'e duality of \cite{CrMo} (Proposition \ref{pdutr}).
This new proof of Poincar\'e duality for leaf spaces appears as a straightforward extension 
of the standard arguments \cite{BoTu} from manifolds to leaf spaces. 
Moreover, this duality extends the known one for basic cohomology of Riemannian foliations \cite{Ser}.\\
\hspace*{.3in}There are several other cohomology theories associated to foliations which are easier to describe and are perhaps more familiar, 
such as basic cohomology (see e.g. \cite{minimal, Ser}) and foliated cohomology (see e.g. \cite{Alv, Hei, KT, MoSo}). In the last two sections 
of our paper, we use our $\check{C}$ech-De Rham model to explicitly describe the relations between the cohomology of leaf spaces
and the basic and foliated cohomology.

\section{Transverse structures on foliations}
\label{transverse}
%***************************************************************
%***************************************************************
%***************************************************************
%***************************************************************

In this section we recall some basic notions concerning the transverse structures of foliations, which formalize
the idea of structures over the leaf space. Throughout, we will work in the smooth context.

\begin{num}{\bf Holonomy}\emph{ Let $M$ be a manifold of dimension $n$, equipped with a foliation $\F$ of codimension $q$.
A {\it transversal section} of $\F$ is an embedded $q$-dimensional submanifold $U\subset M$ which is everywhere transverse to the
leaves. Recall that if $\alpha$ is a path between two points
$x$ and $y$ on the same leaf, and if $U$ and $V$ are transversal sections through $x$ and $y$, then $\alpha$ defines a transport
along the leaves from a neighborhood of $x$ in $U$ to a neighborhood of $y$ in $V$, hence a germ of a diffeomorphism $hol(\alpha): (U, x)\rmap (V, y)$,
called the {\it holonomy} of the path $\alpha$. Two homotopic paths always define the same holonomy. The familiar {\it holonomy groupoid}
\cite{CoOp, Haefl, Wi} is the groupoid $Hol(M, \F)$ over $M$ where arrows $x\rmap y$ are such germs $hol(\alpha)$. If the above transport
``along $\alpha$'' is defined in all of $U$ and embeds $U$ into $V$, this embedding $h: U\hookrightarrow V$ is sometimes also denoted by
$hol(\alpha): U\hookrightarrow V$. Embeddings of this form will be called {\it holonomy embeddings}. Note that composition of paths also induces an 
operation of composition on those holonomy embeddings. (In section \ref{etale} below we will present a more general definition of the so-called
``embedding category'').}
\end{num}

\begin{num}{\bf Transversal basis}\label{rmks2.1}\emph{ Transversal sections $U$ through $x$ as above should be thought of as neighborhoods of the leaf through $x$ in
the leaf space. This motivates the definition of a {\it transversal basis} for $(M, \F)$ as a family $\U$ of transversal sections $U\subset M$ 
with the property that, if $V$ is any transversal section through a given point $y\in M$, there exists a holonomy embedding $h: U\hookrightarrow V$ with $U\in \U$ and $y\in h(U)$.\\
\hspace*{.3in}Typically, a transversal section is a $q$-disk given by a chart for the foliation. Accordingly, we can construct 
a transversal basis $\U$ out of a basis $\tilde{\U}$ of $M$ by domains of foliation charts $\phi_{U}: \tilde{U}\tilde{\rmap} \mathbb{R}^{n-q}\times U$, $\tilde{U}\in \tilde{\U}$, with $U=\mathbb{R}^q$. Note that each inclusion $\tilde{U}\hookrightarrow \tilde{V}$
between opens of $\tilde{\U}$ induces a holonomy embedding $h_{U, V}: U\rmap V$ defined by the condition that the plaque in $\tilde{U}$ through 
$x$ is contained in the plaque in $\tilde{V}$ through $h_{U, V}(x)$.}
\end{num}

\begin{num}\label{trbd}{\bf Transversal bundles}\emph{ Let $G$ be a Lie group and let $\pi: P\rmap M$ be a principal $G$-bundle over $M$. Recall \cite{KaTo}
that $P$ is said to be {\it foliated} if $P$ is equipped with a $G$-equivariant foliation $\tilde{\F}$, of the same dimension as $\F$, whose leaves are transversal to the fibers of $\pi$
and mapped by $\pi$  to those of $\F$. The vectors tangent to $\tilde{\F}$ define a flat partial connection on $P$. In particular, any path $\alpha$
in a leaf $L$ from $x$ to $y$ defines a parallel transport $P_x\rmap P_y$ which depends only on the homotopy class of $\alpha$. We call $P$ a {\it transversal}
principal bundle if the transport depends just on the holonomy of $\alpha$. A vector bundle $E$ on $M$ is said to be foliated (transversal) if the associated principal $GL_r$-bundle is foliated (transversal). By the usual relation between Cartan-Ehresmann connections and Koszul connections, we see that a foliated vector bundle is a vector bundle $E$ over $M$ endowed with a ``flat $\F$-connection'', i.e. an operator
\[ \nabla: \Gamma(\F)\times \Gamma(E)\rmap \Gamma(E) \]
satisfying the usual relations $\nabla_{fX}(s)= f\nabla_{X}(s)$, $\nabla_{X}(fs)= f\nabla_{X}(s)+ X(f)s$, as well as the flatness
relation $\nabla_{[X, Y]}= [\nabla_{X}, \nabla_{Y}]$, for all $X, Y\in \Gamma(\F)$, $f\in C^{\infty}(M)$, $s\in \Gamma(s)$. \\
\hspace*{.3in} Notice that if $P$ is a transversal (principal or vector) bundle, any holonomy embedding $h: U\hookrightarrow V$
induces a well-defined map $h_*: P|_{U} \rmap P|_{V}$, which is functorial in $h$. We will usually just write $h: P|_{U} \rmap P|_{V}$ again for this map.\\
\hspace*{.3in}The basic example of a transversal vector bundle is the normal bundle of the foliation, $\nu= TM/\F$. The associated Koszul connection is precisely the Bott connection \cite{Bott}, $\nabla_{X}(\overline{Y})= \overline{[X, Y]}$. It is a transversal bundle by the very definition of (linear) holonomy.
}
\end{num}

\begin{num}\label{trsh}{\bf Transversal sheaves}\emph{ Analogous definitions apply to sheaves. A sheaf $\A$ on $M$ is called {\it foliated} if its restriction to 
each leaf is locally constant. Thus, (the homotopy class of) a path $\alpha$ from $x$ to $y$ in a leaf $L$ induces an isomorphism between stalks
$\alpha_*: \A_x\rmap \A_y$. The sheaf is transversal if this isomorphism only depends on the holonomy of $\alpha$. A global section
$s\in \Gamma(M, \A)$ is called {\it invariant} if $s$ is invariant under transport along leaves, i.e. $\alpha_*s(x)= s(y)$  in the notations above.\\
\hspace*{.3in}If $\A$ is a transversal sheaf, any holonomy embedding $h: U\hookrightarrow V$ gives a well-defined restriction $h^*: \Gamma(V, \A)\rmap \Gamma(U, \A)$. The global section $s$ is invariant if and only if $h^*(s|_{V})= s|_{U}$ for each such $h$. \\
\hspace*{.3in}An example of a transversal sheaf is the sheaf $\Omega_{bas}^{0}$ of smooth functions which are locally constant along the leaves.
One similarly has the transversal sheaves $\Omega^{k}_{bas}$ of germs of basic differential $k$-forms.
More generally, any foliated vector bundle $E$ induces a foliated sheaf $\Gamma_{\nabla}(E)$ defined as follows.
We denote by $\Gamma_{\nabla}(M; E)$ the space of sections of $E$ which are $\nabla$-constant. Over $M$, $\Gamma_{\nabla}(E)$
is the sheaf whose space of sections over an open $U$ is $\Gamma_{\nabla}(E|_{\, U})$. Using the parallel transport with respect to $\nabla$ 
we see that this sheaf is locally constant when restricted to leaves, hence it is foliated. Clearly $\Gamma_{\nabla}(E)$ is
transversal if $E$ is. For instance, if $E= \Lambda^{k}\nu^*$, then 
$\Gamma_{\nabla}(\Lambda^{k}\nu^*)= \Omega^{k}_{bas}$. \\
\hspace*{.3in}Another important example is the (real) transversal orientation sheaf of the foliation, which we denote by $\oo$.
When restricted to a transversal open $U$, $\Gamma(U; \oo)= H^{q}_{c}(U)^{\vee}$. The foliation is transversally
orientable if and only if $\oo$ is constant. 
}
\end{num}

%***************************************************************
%***************************************************************
%***************************************************************
%***************************************************************
\section{The transversal $\check{C}$ech-De Rham complex}
\label{CDRcomplex}
%***************************************************************
%***************************************************************
%***************************************************************
%***************************************************************

Let $(M, \F)$ be a foliated manifold and let $\U$ be a transversal basis. Consider the double complex
which in bi-degree $k, l$ is the vector space 
\[ C^{k, l}= \check{C}^{k}(\U, \Omega^{l})= \prod_{U_{0}\stackrel{h_1}{\rmap} \ldots \stackrel{h_k}{\rmap} U_k} \Omega^l(U_0) .\]
Here the product ranges over all $k$-tuples of holonomy embeddings between transversal sections from the given basis $\U$,
and $\Omega^k(U_0)$ is the space of differential $k$-forms on $U_0$. For elements $\omega\in C^{k, l}$, we denote its components by
$\omega(h_1, \ldots, h_k)\in \Omega^k(\U_0)$. The vertical differential $C^{k, l}\rmap C^{k, l+1}$ is $(-1)^kd$ where $d$ is the usual De Rham differential.
The horizontal differential $C^{k, l}\rmap C^{k+1, l}$ is $\delta= \sum(-1)^{i}\delta_{i}$ where
\begin{equation}\label{deltas} \delta_{i}(h_1, \ldots , h_{k+1})= \left\{ \begin{array}{lll}
                                      h_{1}^{*}\omega(h_2, \ldots , h_{k+1}) \ \ \mbox{if $i=0$}\\ 
                                      \omega(h_1, \ldots, h_{i+1}h_{i}, \ldots, h_{k+1}) \ \ \mbox{if $0<i< k+1$}\\
                                      \omega(h_1, \ldots, h_k) \ \ \mbox{if $i= k+1$}
                        \end{array}
                \right.
\end{equation}
This double complex is actually a bigraded differential algebra, with the usual product 
\[ (\omega\cdot\eta)(h_1, \ldots , h_{k+k\,'})= (-1)^{kk\,'}\omega(h_1, \ldots , h_{k}) h_{1}^{*} \ldots h_{k}^{*} \eta(h_{k+1}, \ldots h_{k+k\,'}) \]
for $\omega\in C^{k, l}$ and $\eta\in C^{k\,', l\,'}$. We will also write $\check{C}(\U, \Omega)$
for the associated total complex, and refer to it as the {\it $\check{C}$ech-De Rham complex} of the foliation. The associated cohomology is denoted 
\[ \check{H}_{\U}^{*}(M/\F) ,\]
and referred to as the $\check{C}$ech-De Rham cohomology of the leaf space $M/\F$, w.r.t. the cover $\U$. \\
\hspace*{.3in}Note that, when $\F$ is the codimension $n$ foliation by points, then $\U$ is a basis for the topology of $M$, and
$C^{k, l}$ is the usual $\check{C}$ech-De Rham complex \cite{BoTu}. Thus in this case $\check{H}_{\U}^*(M/\F)= H^*(M)$ is the usual De Rham cohomology
of $M$.\\
\hspace*{.3in}In general, choosing a transversal basis $\U$ and a basis $\tilde{\U}$ of $M$ as in \ref{rmks2.1}, there is
an obvious map of double complexes $C^{k, l}(\U) \rmap C^{k, l}(\tilde{\U})$ into the $\check{C}$ech-De Rham complex for the manifold $M$. Hence a canonical map
\begin{equation}\label{zero}\label{map2.1}
\pi^{*}: \check{H}_{\U}^{*}(M/\F)\rmap H^{*}(M; \mathbb{R}) \ ,
\end{equation}
which should be thought of as the pull-back along the ``quotient map'' $\pi: M\rmap M/\F$.\\

\hspace*{.1in}The standard way \cite{Co3, minimal} to model the leaf space of a foliation $(M, \F)$ is by the classifying space $BHol(M, \F)$ of the holonomy groupoid. Thus, the
following theorem can be interpreted as a $\check{C}$ech-De Rham theorem for leaf spaces.

\begin{st}\label{theorem1}There is a natural isomorphism 
\[ \check{H}_{\U}^{*}(M/\F) \cong H^{*}(B Hol(M, \F); \mathbb{R})\ ,\]
between the $\check{C}$ech-De Rham cohomology and the cohomology of the classifying space. In particular, the left hand side is independent
of the choice of a transversal basis $\U$.
\end{st}

\hspace*{.1in}For the proof of this theorem, we choose a complete transversal section $T$ which contains every $U\in \U$,
and we consider the ``reduced holonomy groupoid'' $Hol_{T}(M, \F)$, defined as the restriction of $Hol(M, \F)$ to $T$. 
We may assume that $\U$ is a basis for the topology of $T$. By a well known Morita-invariance argument, the classifying spaces $B Hol(M, \F)$ and $B Hol_{T}(M, \F)$
are weakly homotopyc equivalent. The advantage of passing to a complete transversal is that $Hol_{T}(M, \F)$ becomes an \'etale groupoid
(see section \ref{etale} for the precise definitions).
For such groupoids $\G$ there is a standard cohomology $H^{*}(\G; -)$ with coefficients, which was also defined by Haefliger \cite{difcoh} 
in terms of bar-complexes, and which is known \cite{conj} to be isomorphic to the cohomology of the classifying space. 
In section \ref{etale} we will recall all the basic definitions. The theorem will then follow from the following proposition,   
which is a particular case of the Theorem \ref{theoreml} below.

\begin{prop}\label{lema1} For any complete transversal $T$ and any basis $\U$ of $T$, there is a natural isomorphism
\[  \check{H}_{\U}^{*}(M/\F) \cong H^{*}(Hol_{T}(M, \F); \mathbb{R}) \ .\]
\end{prop}

\hspace*{.1in}We mention here that there are several variations of Theorem \ref{theorem1}. For instance, for any 
 transversal sheaf $\A$ there is a $\check{C}$ech complex $\check{C}(\U, \A)$. In degree $k$,
\[ \check{C}^{k}(\U; \A)= \prod_{U_{0}\stackrel{h_1}{\rmap} \ldots \stackrel{h_k}{\rmap} U_k} \Gamma(U_0; \A) \ ,\]
with the boundary $\delta= \sum(-1)^{i}\delta_{i}$ given by the formulas (\ref{deltas}). A consequence of the more general Theorem \ref{theoreml} says that, if $\A|_{U}$ is acyclic for all $U\in \U$, then $\check{C}(\U, \A)$ computes the cohomology of the classifying space (of the reduced holonomy groupoid) with coefficients in a sheaf $\tilde{\A}$ naturally associated to $\A$.\\
\hspace*{.3in}Another variation applies to the cohomology with compact supports (see section \ref{etale}). Note that all these
are actually extensions of the usual ``$\check{C}$ech-De Rham isomorphisms'' \cite{BoTu} from manifolds to leaf space. Accordingly, an immediate consequence will be the Poincar\'e duality for leaf spaces 
(see Section \ref{secbasic}), which is one of the
main results of \cite{CrMo}. With Theorem \ref{theorem1} and its analogue for compact supports available, the new proof of Poincar\'e duality is this time a rather straightforward extension of the classical proof
\cite{BoTu} from manifolds to leaf spaces.

%***************************************************************
%***************************************************************
%***************************************************************
%***************************************************************
\section{The transversal Chern-Weil map}
\label{classes}
%***************************************************************
%***************************************************************
%***************************************************************
%***************************************************************

To illustrate the usefulness of the transversal $\check{C}$ech-De Rham complex we will adapt the standard geometric construction
of characteristic classes of principal bundles to transversal bundles, so as to obtain explicit classes in this complex. We will use the Weil-complex formulation, which we recall first (for an extensive exposition,
see \cite{KaTo, Duff}).

\begin{num}\label{clCW}{\bf Classical Chern-Weil: }\rm \cite{Car} Recall that the Weil algebra of the Lie algebra $\mathfrak{g}$ (of a Lie group $G$) is the algebra \[ W(\mathfrak{g})= S(\mathfrak{g}^*)\otimes \Lambda(\mathfrak{g}^*) .\]
It is a graded commutative dga (graded as $W(\mathfrak{g})^n= \oplus_{2p+q= n} S^p(\mathfrak{g}^*)\otimes \Lambda^q(\mathfrak{g}^*))$, equipped
with operations $i_X$ and $L_X$ (linear in $X\in \mathfrak{g}$) which satisfy the usual Cartan identities. In the language of \cite{KaTo}, 
this means that $W(\mathfrak{g})$ is a $\mathfrak{g}$-dga. If $P$ is a principal $G$-bundle over a manifold $M$, the algebra $\Omega^*(P)$ of differential forms on $P$ with its usual operations $i_X$ and $L_X$ is another example of a $\mathfrak{g}$-dga. A connection $\nabla$ on $P$ is uniquely determined by its connection form $\omega\in \Omega^1(P)\otimes\mathfrak{g}$. This can be viewed as a map $\omega: W(\mathfrak{g})^1= \mathfrak{g}^*\rmap \Omega^1(P)$, which extends uniquely to a map of $\mathfrak{g}$-dga's, 
\begin{equation}
\label{charact}
\tilde{k}(\nabla): W(\mathfrak{g}) \rmap \Omega(P) \ .
\end{equation}
(On $\mathfrak{g}^*= S^1(\mathfrak{g}^*)\subset W(\mathfrak{g})^2$, it restricts to the curvature $\Omega= d\omega+ \frac{1}{2} [\omega, \omega]$.) The restriction of this map (\ref{charact}) to basic elements (elements annihilated by $i_X$ and $G$-invariant) gives a map of dga's
\[ S(\mathfrak{g}^*)^{G} \rmap \Omega^*(M)\]
(zero differential on $S(\mathfrak{g}^*)^{G}$, the usual De Rham differential on $\Omega(M)$), hence a map
\begin{equation}\label{classicalCW} 
k(\nabla): S(\mathfrak{g}^*)^{G} \rmap H^{*}(M) ,
\end{equation}
known as the Chern-Weil map for the principal $G$-bundle $P$. Because of the $2p$ in the
grading of the Weil algebra, $k(\nabla)$ maps invariant polynomials of degree $p$ to degree $2p$ cohomology classes. Moreover, $k(\nabla)$
does not depend on $\nabla$. This follows from the Chern-Simons construction (see below). A more refined characteristic map is obtained 
if one uses a maximal compact subgroup $K$ of $G$. Since $P/K\rmap M$ has contractible fibers, the map induced in De Rham cohomology is
an isomorphism. Hence, to get down to the base manifold, it suffices to consider the $K$-basic elements of (\ref{charact}). Denoting
by $W(\mathfrak{g}, K)$ the subcomplex of $W(\mathfrak{g})$ of $K$-basic elements, one obtains a characteristic map $H^*(W(\mathfrak{g}, K))\rmap H^{*}(M)$.
\end{num}

\begin{num}\label{Simons}{\bf Chern-Simons: }\rm Given $k$ connections $\nabla_0, \ldots , \nabla_k$ on $P$, we consider the convex combination 
\begin{equation}
\label{doi}
\nabla= t_0\nabla_0+ \ldots + t_k\nabla_k
\end{equation}
which defines a connection on the principal bundle ${\bf \Delta}^{k}\times P$ over 
${\bf \Delta}^{k}\times M$, where ${\bf \Delta}^{k}= \{ (t_0, \ldots, t_k): t_i\geq 0, \sum t_{i}= 1\}$ is the standard $k$-simplex.
We define
\begin{equation}\label{patruu}
\tilde{k}(\nabla_0, \ldots, \nabla_k)= (-1)^{k}\int_{{\bf \Delta}^{k}}\tilde{k}(\nabla): W(\mathfrak{g}) \rmap \Omega^{*-k}(P)\ ,
\end{equation}
where $\int_{{\bf \Delta}^{k}}: \Omega^*({\bf \Delta}^{k}\times P)\rmap \Omega^{*-k}(P)$ is the integration along the fibers ${\bf \Delta}^{k}$.
Let us summarize the main properties of this construction:
\begin{enumerate}[(i)]
\item the map (\ref{patruu}) commutes with the action of $G$, and with the operators $i_X$, $L_X$, and it vanishes on all elements 
$\alpha\otimes \beta\in W(\mathfrak{g})$ with $\alpha$ a polynomial of degree $> dim(M)$.
\item \begin{equation}\label{Stokes} 
[\tilde{k}(\nabla_0, \ldots, \nabla_k), d]= \sum_{i=0}^{k} (-1)^{i} \tilde{k}(\nabla_0, \ldots , \widehat{\nabla_i}, \ldots , \nabla_k)\ ,
\end{equation} 
\item (\ref{patruu}) is natural w.r.t. isomorphisms of principal $G$-bundles.
\end{enumerate}
\end{num}

%\begin{lem}\label{lemSimons}\ \ 
%\begin{enumerate}[(i)]
%\item the map (\ref{patru}) commutes with the action of $G$, and with the operators $i_X$, $L_X$, and it vanishes on all elements 
%$\alpha\otimes \beta\in W(\mathfrak{g})$ with $\alpha$ a polynomial of degree $> dim(M)$.
%\item \begin{equation}\label{Stokes} 
%[\tilde{k}(\nabla_0, \ldots, \nabla_k), d]= \sum_{i=0}^{k} (-1)^{i} \tilde{k}(\nabla_0, \ldots , \widehat{\nabla_i}, \ldots , \nabla_k)\ ,
%\end{equation} 
%\item $\tilde{k}(h^*\nabla_0, h^*\ldots, h^*\nabla_k)= h^*\tilde{k}(\nabla_0, \ldots, \nabla_k)$, for any isomorphism $h: P\rmap Q$ of principal %$G$-bundles.
%\end{enumerate}
%\end{lem}

{\it Proof:} (ii) is just a version of Stokes' formula (see also \cite{Bott}), while (iii) is obvious. We prove the vanishing result of (i).
Denote by $d$ the degree of the polynomial $\alpha$ and by $q$ the dimension of $M$. We prove that when $d< k$
or $2d> q+k$, our expression
\[ \theta= \tilde{k}(\nabla_0, \ldots, \nabla_k)(\alpha\otimes\beta) \]
vanishes (note that if $d> q$, then at least one of these two equalities holds). First assume that $d<k$. We have $\tilde{k}(\nabla)(\alpha\otimes \beta)= \alpha(\Omega)\wedge \beta(\omega)$, where $\nabla$ is the affine 
combination (\ref{doi}), $\omega$ is the associated $1$-form, and $\Omega$ is its curvature. Let us say that a homogeneous form 
$f dt_{i_1}\ldots dt_{i_r}dx_{j_1}\ldots dx_{j_s}$ on ${\bf \Delta}^{k}\times P$ has bi-degree $(r, s)$. Since $\omega$ has bi-degree $(0, 1)$, 
$\Omega$ is a sum of forms
of bi-degree $(1, 1)$ and $(0, 2)$, so $\int_{{\bf \Delta}^{k}} \alpha(\Omega)\wedge \beta(\omega)= 0$ because no bi-degree
$(r, s)$ with $r= k$ will be involved.\\
We now turn to the case $2d> q+k$. Let $l$ be the degree of $\beta$. Because of the similar property for $\beta$, we have
$i_{X_1}\ldots i_{X_{l+1}}\theta = 0$ for any vertical vector fields $X_i$. On the other hand, $i_{Y_1}\ldots i_{Y_{q+1}}\theta = 0$ for
any horizontal vector fields $Y_i$. Since $deg(\theta)= 2d+l-k> l+ q$, it follows that $\theta= 0$. \ \ $\Boxe$

\begin{num}\label{trChW}{\bf Construction of the transversal Chern-Weil map: }\rm Now let $P$ be a transversal principal $G$-bundle on a foliated manifold $(M, \F)$. Consider the $\check{C}$ech-De Rham complex 
\[ \check{C}^{k}(\U, \Omega^{l}(P))= \prod_{U_{0}\stackrel{h_1}{\rmap} \ldots \stackrel{h_k}{\rmap} U_k} \Omega^l(P|_{U_{0}})\ ,\]
defined exactly as in section \ref{CDRcomplex} (except that $\Omega^l(U_0)$ is replaced by $\Omega^l(P|_{U_{0}})$, and hence the horizontal differential $\delta$ involves the maps $h_{1}: P|_{U_{0}}\rmap P|_{U_{1}}$ discussed in section \ref{transverse}). Choose a system $\nabla= \{ \nabla_{U}\}$ of connections, one connection $\nabla_{U}$ on $P|_{U}$ for each $U$ in a transversal basis $\U$. 
In general we cannot assume this choice to be respected by holonomy embeddings $h: U\rmap V$, i.e. $\nabla_U$ is in general different from the connection on $P|_{U}$ induced by $\nabla_{V}$ via the isomorphism $h: P|_{U} \rmap P|_{h(U)}$. Denote this last connection by $\nabla_{h}$. For a string $U_{0}\stackrel{h_1}{\rmap} \ldots \stackrel{h_k}{\rmap} U_k$ of holonomy embeddings, we consider the map
(see \ref{clCW} above) 
\begin{equation}\label{patru}
\tilde{k}(\nabla_{U_0},
\nabla_{h_1}, \nabla_{h_2h_1}, \ldots ,
\nabla_{h_k \ldots h_2h_1}): W(\mathfrak{g}) \rmap \Omega^{*-k}(P|_{U_{0}}) .
\end{equation}
Doing this 
for all such strings, we obtain a map into the total complex 
\begin{equation}\label{cinci}
\tilde{k}(\omega): W(\mathfrak{g}) \rmap \prod_{U_{0}\stackrel{h_1}{\rmap} \ldots \stackrel{h_k}{\rmap} U_k} \Omega^{*-k}(P|_{U_{0}})\ .
\end{equation}
This map respects the total degree, and it is obviously compatible with the operations $i_X$ and the $G$-action. So, by restricting to basic elements
it yields a map into the transversal $\check{C}$ech-De Rham complex
\begin{equation}\label{sase}
\tilde{k}(\omega): S(\mathfrak{g}^{*})^{G} \rmap \check{C}^{*}(\U, \Omega^*) 
\end{equation}
(mapping degree $p$ polynomials into elements of total degree $2p$). 
\end{num}

\begin{st}\label{theorem2} The Chern-Weil map of a transversal principal $G$-bundle $P$ over $(M, \F)$ 
has the following properties:
\begin{enumerate}[(i)]
\item The maps {\rm (\ref{cinci})} and {\rm (\ref{sase})} respect the differential, hence they induce a map
\begin{equation}\label{sapte}
k_{P}:= k(\nabla) : S(\mathfrak{g}^{*})^{G} \rmap \check{H}^{*}_{\U}(M/\F)\ ,
\end{equation}
\item This map {\rm (\ref{sapte})} does not depend on the choice of the connections $\{\nabla_{U}\}$, and respects the products.
\item Composed with the pull-back map $\pi^{*}: \check{H}^{*}_{\U}(M/\F)\rmap H^{*}(M)$, see {\rm (\ref{map2.1})}, it gives the usual Chern-Weil map {\rm (\ref{classicalCW})} of $P$. 
\item {\rm (}``Bott vanishing theorem''{\rm )} The image of the map {\rm (\ref{sapte})} 
is zero in degrees $> 2q$, where $q$ is the codimension of $\F$.
\end{enumerate}
\end{st}

\hspace*{.1in}The classical Bott vanishing theorem \cite{Bott} (for the normal bundle of the foliation) and its extensions
to foliated bundles \cite{KaTo} are at the level of $H^*(M)$. The point of Theorem \ref{theorem2} is that,
using {\it classical geometrical arguments}, one can prove these vanishing results and construct the resulting
characteristic classes at the level of the leaf space, i.e. in the cohomology of the classifying space (cf. Theorem \ref{theorem1}).\\

%Of course, this only applies to the foliated bundles which are transversal (the reason being that only these bundles come from
%the holonomy groupoid of $\F$). In the next section we will extend our arguments to general \'etale groupoids. Hence a consequence
%of this is that Theorem \ref{theorem2} and Corollary \ref{corex} does hold for general foliated bundles, provided
%we replace  $\check{H}_{\U}^{*}(M/\F)\cong H^*(B Hol(M, \F))$ by the cohomology of the (classifying space of the) monodromy groupoid of $\F$.

{\it Proof of Theorem \ref{theorem2}:} (i) and (iv) clearly follow from the main properties of the Chern-Simons construction \ref{Simons}. Also (iii) will follow from the independence of the connections. Indeed, it suffices to check that, if $\F$ is the foliation by points, then the resulting map $k_{\nabla}: S(\mathfrak{g}^{*})^{G}\rmap 
\check{H}^{*}_{\U}(M)$ composed with $\check{C}$ech-De Rham isomorphism $\check{H}^{*}_{\U}(M)\cong H^{*}(M)$ (induced by the inclusion
$\Omega^*(M)\subset C^*(\U, \Omega^*)$ \cite{BoTu}) gives the usual Chern-Weil map. But this is clear even at the chain level, provided we choose 
$\nabla_{U}= \nabla |_{U}$ for some globally defined connection $\nabla$.\\
(ii) For two different choices $\nabla= \{\nabla_{U}\}$ and $\nabla\,'= \{\nabla_{U}\,'\}$ of connections, the map
$H: W(\mathfrak{g})\rmap C^*(\U,\Omega^*)$ defined by
\[ H^{*}(w)(h_1, \ .\ .\ .\ , h_k)= 
\sum_{i=0}^{k}(-1)^i k(\nabla_{h_i \ldots h_2h_1}, \ldots , \nabla_{h_k \ldots h_2h_1}, \nabla^{\,'}_{U_0}, \nabla^{\,'}_{h_1},  \ldots , 
\nabla^{\, '}_{h_i \ldots h_2h_1})(w) \ . \]
provides an explicit chain homotopy. To prove the compatibility with the products, one can either proceed as in \cite{KaTo} using the simplicial
Weil complex (see \cite{Crath} for details), or, since the characteristic map is constructed through the double complex 
$\check{C}^{p}(\U, \Omega^{p+q}({\bf \Delta}^{q} \times P))$ by integration over the simplices, one can use the simplicial De Rham complex and Theorem 2.14 of \cite{Dupon}.  \ \ $\Boxe$

\begin{num}\label{exotic}{\bf Exotic characteristic classes: }\rm The vanishing result of Theorem \ref{theorem2} shows that the
construction of the ``exotic'' classes also lifts to the $\check{C}$ech-De Rham complex. To describe all the relevant characeristic classes,
we consider the complex $W(\mathfrak{g}, K)$ of $K$-basic elements described in \ref{clCW}, together with its $q$-th truncation 
$W_q(\mathfrak{g}, K)$ defined as the quotient by the ideal generated by the elements of polynomial degree $> q$. 
By the vanishing result (more precisely from the proof above), 
the map (\ref{cinci}) induces a chain map $W_q(\mathfrak{g}, K) \rmap \check{C}^{*}(\U, \Omega^*(P/K)$. 
Using the contractibility of $G/K$ as in 
\ref{clCW}, we obtain the following refinement of the characteristic map of Theorem \ref{theorem2}.
\end{num}

\begin{cor}\label{corex} The Chern-Weil construction of {\rm \ref{trChW}} gives a well-defined algebra map
\begin{equation}\label{exoticone}
k^{ex}_{P}:= k^{ex}(\nabla): H^*(W_q(\mathfrak{g}, K)) \rmap \check{H}_{\U}^{*}(M/\F) \ ,
\end{equation}
again independent of the choice of connections. Moreover, composed with the pull-back map $\pi^{*}: \check{H}^{*}_{\U}(M/\F)\rmap H^{*}(M)$ (see {\rm (\ref{map2.1})}), it gives the exotic characteristic map of the foliated bundle $P$ {\rm \cite{KaTo}}.
\end{cor}

%***************************************************************
%***************************************************************
%***************************************************************
%***************************************************************
\section{The $\check{C}$ech-De Rham complex of an \'etale groupoid}
\label{etale}
%***************************************************************
%***************************************************************
%***************************************************************
%***************************************************************

In this section we prove Theorem \ref{theorem1}, as well as some generalizations and variants, in the context of \'etale groupoids.
Our general goal is to describe the (hyper-) homology and cohomology of \'etale groupoids in terms of the (hyper-) homology and
cohomology of small categories. We begin by introducing some standard terminology.

\begin{num}{\bf Smooth \'etale groupoids: }\rm A {\it smooth groupoid} is a groupoid $\G$ for which the sets $\nG{0}$ and $\nG{1}$ of objects 
and arrows have the structure of a smooth manifold, all the structure maps are smooth, and the source and the target maps are moreover submersions.
The holonomy groupoid $Hol(M, \F)$ of a foliation is an example of a smooth groupoid. Such a smooth groupoid is said to be \'etale if the source
and the target maps are local diffeomorphisms. In this case the manifolds $\nG{0}$ and $\nG{1}$ have the same dimension, to which we refer as 
{\it the dimension of} $\G$. An example of an \'etale groupoid of dimension $q$ is the universal Haefliger groupoid $\Gamma^q$ for codimension $q$ foliations \cite{Haefl}. There is an important notion of {\it Morita equivalence} between smooth groupoids, see e.g. \cite{Co3, Haefl, Fourier, Mrcun}.
For any foliation, the holonomy groupoid $Hol(M, \F)$ is Morita equivalent to an \'etale groupoid, namely to its restriction to any
complete transversal $T$, denoted $Hol_{T}(M, \F)$. A Morita equivalence between smooth groupoids induces a weak homotopy equivalence between their classifying spaces.
\end{num}

\begin{num}\label{bar}{\bf Sheaves and cohomology: }\rm For a smooth \'etale groupoid $\G$, a $\G$-sheaf is a sheaf $\A$ over the space $\nG{0}$, equipped with a continuous $\G$-action. For any such sheaf there are natural cohomology groups $H^n(\G; \A)$ whose definition we recall.
Denote by $\nG{k}$ the space of composable arrows
\begin{equation}\label{string}
x_0\stackrel{g_1}{\rmap} \ldots \stackrel{g_k}{\rmap} x_k
\end{equation}
of $\G$, and by $\epsilon_{k}: \nG{k}\rmap \nG{0}$ the map which sends (\ref{string}) 
to $x_0$. The bar complex of $\A$ is defined by $B^k(\G; \A)= \Gamma(\nG{k}; \epsilon_{k}^{*}\A)$, hence consists on continuous functions $c$
which associate to a string of arrows (\ref{string}) an element $c(g_1, \ldots, g_k)\in \A_{x_0}$. The boundary is $\delta= \sum (-1)^{i}\delta_i$ 
with the same formulas as in (\ref{deltas}). If $\A$ is ``good'' in the sense that $\A$ and its pull-backs $\epsilon_{k}^{*}\A$ are injective sheaves, then $H^n(\G; \A)$ is computed by the bar complex $B(\G, \A)$. In general, one chooses a resolution $\es^{*}$ of $\A$ by ``good'' $\G$-sheaves,
and $H^n(\G; \A)$ is computed by the double complex $B^k(\G; \es^l)$. In general, these cohomology groups coincide with the cohomology groups of the classifying space $B\G$ \cite{conj}. \\
\hspace*{.3in}Similarly, using compact supports and direct sums in the definition of the bar complex, one defines the homology groups $H_{*}(\G; \A)$  \cite{CrMo} (sometimes denoted $H_{c}^{*}(\G; \A)= H_{-*}(\G; \A)$), which should be thought of as
a good model for the compactly supported cohomology of the classifying space.   
\end{num}

\begin{num}{\bf $\check{C}$ech complexes: }\rm Let $\G$ be an \'etale groupoid and let $\U$ be a basis of opens in $\nG{0}$. 
A {\it $\G$-embedding} $\sigma: U\rmap V$ is a smooth family $\sigma(x)$, $x\in U$, where each $\sigma(x): x\rmap y$ is an arrow in $\G$ from $x$ to some point $y\in V$; moreover, the map $x\rmap $target$(\sigma(x))$ should define an embedding of $U$ into $V$. As in the first section, we can
now define the $\check{C}$ech complex $\check{C}(\U; \A)$ for any $\G$-sheaf $\A$, 
\[ \check{C}^{k}(\U; \A)= \prod_{U_0 \rmap \ldots \rmap U_k} \Gamma(U_0, \A)\ ,\]
where the product is over all strings of $\G$-embeddings between opens $U\in \U$, 
and the boundary $\delta= \sum(-1)^{i} \delta_i$ is given by the same formulas as in (\ref{deltas}). \\
\hspace*{.3in}We say that $\A$ is $\U$-acyclic if $H^{i}(U; \A)= 0$ for each $i>0$ and each $U\in \U$. In this case
define $\check{H}_{\U}^{*}(\G; \A)$ as the cohomology of $\check{C}(\U; \A)$. In general, we define
$\check{H}_{\U}^{*}(\G; \A)$ as the cohomology of the double complex $\check{C}^{k}(\U; \es^{l})$, where 
$0\rmap \A\rmap \es^0\rmap\ldots \rmap \es^d\rmap 0$ is a bounded resolution by $\U$-acyclic sheaves, 
$d= dim(\G)$. By the usual arguments, such resolutions always exist, and the definition does not depend on the choice of
the resolution.
\end{num}

\begin{num}\label{CDRetale}{\bf Examples: }\rm  The $\G$-sheaf $\Omega_{\G}^{l}$ of $l$-differential forms with its natural $\G$-action is always $\U$-acyclic, as is any soft $\G$-sheaf. We obtain the $\check{C}$ech-De Rham (double) complex of $\G$, $\check{C}(\U; \Omega)$, computing $\check{H}^{*}_{\U}(\G; \mathbb{R})$. 
If the basis $\U$ consists of contractible opens (balls), then any locally constant $\G$-sheaf $\A$ is $\U$-acyclic, hence
$\check{H}^{*}_{\U}(\G; \A)$ is computed by $\check{C}(\U; \A)$.
\end{num}

\hspace*{-.2in}Similarly one defines the $\check{C}$ech complex with compact supports $\check{C}_{c}^{*}(\U; \A)$ using
\[ \bigoplus_{U_0 \rmap \ldots \rmap U_k} \Gamma_{c}(U_0, \A)\ .\]
In order to get a cochain complex, we associate the degree $-k$ to the direct sums over strings of $k$ $\G$-embeddings. 
If $\A$ is $c$-soft, then $\check{H}_{c,\, \U}^{*}(\G; \A)$ is defined by $\check{C}_{c}(\U; \A)$. In general,
one uses a resolution  $0\rmap \A\rmap \es^0\rmap\ldots \rmap \es^d\rmap 0$ by $c$-soft $\G$-sheaves, and the
double complex $\check{C}^{k}_{c}(\U; \es^{l})$. The resulting cohomology is denoted $\check{H}^{*}_{c,\, \U}(\G; \A)$.

\begin{num}{\bf The embedding category: }\label{embdcat}\rm The notion of $\G$-embedding originates in \cite{embd}, where the second author has introduced a small category $Emb_{\U}(\G)$ for each basis $\U$ of open sets. The objects of $Emb_{\U}(\G)$ are the members $U$ of $\U$, and the arrows are the $\G$-embeddings between the opens of $\U$. The main result of \cite{embd} was that the classifying space $B\G$ is 
weakly homotopy equivalent to the CW-complex $BEmb_{\U}(\G)$, provided each of the basic opens in $\U$ is contractible. \\
\hspace*{.3in}Now any $\G$-sheaf $\A$ defines 
an obvious contravariant functor $\Gamma(\A)$ on $Emb_{\U}(\G)$
sending $U$ to $\Gamma(U; \A)$, and $\check{C}(\U; \A)$ is just the usual (bar) complex computing
the cohomology $H^*(Emb_{\U}(\G); \Gamma(\A))$ of the discrete category $Emb_{\U}(\G)$ with coefficients. 
Hence \cite{embd} proves that $H^*(\G; \A)= \check{H}^{*}_{\U}(\G; \A)$ provided all the opens $U\in \U$ are contractible
and $\A$ is (locally) constant. We now prove a stronger ``$\check{C}$ech-De Rham isomorphism''
which applies to more general coefficients, and also to compact supports.
\end{num}

\begin{st}\label{theoreml} Let $\G$ be an \'etale groupoid, and let $\U$ be a basis for $\nG{0}$ as above. Then for any $\G$-sheaf $\A$, there
are natural isomorphisms
\[  H^n(\G; \A)= \check{H}^{n}_{\U}(\G; \A) ,\ \  H^{n}_{c}(\G; \A)= \check{H}^{n}_{c,\, \U}(\G; \A) \ .\] 
\end{st}

{\it Proof:} The proofs of the isomorphisms in the statement are similar, and we only prove the first one (an explicit proof of the second one also occurs in \cite{Crath}). By comparing
resolutions of the $\G$-sheaf $\A$, it suffices to find a suitable complex $C(\A)$ and
explicit quasi-isomorphisms 
\[ B(\G; \A) \lmap C(\A) \rmap \check{C}(\U; \A) \]
natural in $\A$, for the case where $\A$ is ``good'' in the sense of \ref{bar}.
For this we consider the bisimplicial space $S_{p, q}$, whose $p, q$-simplices are of the form
\begin{equation}\label{doublestring}
 x_0\stackrel{g_1}{\rmap} \ldots \stackrel{g_q}{\rmap} x_q\stackrel{g}{\rmap} U_0 \stackrel{\sigma_1}{\rmap} \ldots \stackrel{\sigma_p}{\rmap} U_p \ ,
\end{equation}
where $\sigma_1, \ldots , \sigma_p$ are $\G$-embeddings, and $g_1, \ldots , g_q, g$ are arrows in the groupoid $\G$, the notation $x_q\stackrel{g}{\rmap} U_0$ indicating that the target of $g$ is in $U_0$. The topology on $S_{p, q}$ is the topology induced from
the topology on $\G$, 
\[ S_{p, q}= \coprod_{U_0\rmap \ldots \rmap U_p} \nG{q} \ .\]
The $\G$-sheaf $\A$ induces a sheaf $\A_{p, q}$ on $S_{p, q}$ by pull-back along the projection which maps (\ref{doublestring}) to $x_q$. Consider the
double complex $C= C(\A)$,
\[ C^{p, q}= \Gamma(S_{p, q}, \A_{p, q}) \ .\]
For a fixed $p$, the complex $C^{p, *}$ is a product of complexes, namely, for each string $U_0\rmap \ldots\rmap U_p$, the bar 
complex (see \ref{bar}) of the (\'etale) comma groupoid $\G/U_0$ with coefficients in the pull-back of the sheaf $\A$. Since the groupoid $\G/U_0$ is Morita equivalent to the space $U_0$, this cohomology is $H^{*}(U_0; \A)$. Since $\A$ is assumed to be good, $H^*(U_0;\A)$ vanishes in positive degrees, and we conclude that the canonical map
\[  \check{C}^{p}(\U; \A) \rmap C^{p, *} \]
is a quasi-isomorphism for each fixed $p$. Write $\pi_{p, q}: S_{p, q}\rmap \nG{q}$ for the projection of (\ref{doublestring}) to the string $x_0 \rmap \ldots \rmap x_q$. Then $C^{p, q}= \Gamma((\pi_{p, q})_*(\A_{p, q}))$. The stalk of $(\pi_{p, q})_*(\A_{p, q})$ at $x_0 \rmap \ldots \rmap x_q$
is 
\begin{equation}\label{colimit}
\lim_{\overrightarrow{x_q\in U}}( \prod_{U\rmap U_0\rmap \ldots \rmap U_p} \Gamma(U; \A) )\ ,
\end{equation}
where the colimit is taken over all basic open neighborhoods $U$ of $x_q$. For a fixed $U$, the complex inside the $\lim$ in (\ref{colimit})
computes the cohomology of the (discrete) comma category $U/Emb_{\U}(\G)$ with coefficients in the constant group $\Gamma(U, \A)$. Since the comma category is contractible, so is this complex. Taking the colimit, we see that for each $q$ the map $\A\rmap (\pi_{p, q})_*(\A_{p, q})$ is a quasi-isomorphism of (complexes of)
injective sheaves on $\nG{q}$. Thus the natural map
\[ B^{q}(\G; \A)= \Gamma(\nG{q}; \epsilon_{q}^{*}\A) \rmap \Gamma((\pi_{*, q})_*(\A_{*, q}))= C^{*, q} \]
is a quasi-isomorphism, and the proof is complete.  \ \ $\Boxe$\\

\hspace*{.1in}Regarding the relation with the embedding category \ref{embdcat} and its cohomology, let us point out 
the following immediate consequence, which
is an improvement of the result of \cite{embd}.\\

\begin{cor} If $\G$ is an \'etale groupoid, $\tilde{\U}$ is a basis of opens of $\nG{0}$, and $\A$ is a $\G$-sheaf
with the property that $H^k(U, \A|_{\,U})= 0$ for all $U\in \U$, $k\geq 1$, then
\[ H^*(\G; \A) \cong H^{*}(Emb_{\U}(\G); \Gamma(\A))\ .\]
Similarly, if $H^{k}_{c}(U, \A|_{\,U})= 0$ for all $U\in \U$, $k\geq 1$, then
\[ H^{*}_{c}(\G; \A) \cong H_{-*}(Emb_{\U}(\G); \Gamma_{c}(\A))\ .\]
Also, if each $U\in \U$ is contractible, and $\A$ is locally constant as a sheaf on $\nG{0}$, then
\[ H^*(\G; \A) \cong H^{*}(Emb_{\U}(\G); \Gamma(\A)),\ \ H^{*}_{c}(\G; \A) \cong H_{d-*}(Emb_{\U}(\G); H^{d}_{c}(-; \A)) \]
(where $d$ is the dimension of the base space $\nG{0}$).
\end{cor}

%\hspace*{.1in}An immediate consequence of Theorem \ref{theoreml} is a new proof of the Poincar\'e duality for \'etale groupoids
%(which is one of the main results in \cite{CrMo}), which appears as an extension of the classical proof for manifolds \cite{BoTu}.
%Let us denote by $\oo$ the (real) orientation sheaf of $\nG{0}$, viewed as a $\G$-sheaf in the natural way. 

%\begin{cor}\label{poincare} (Poincare duality \cite{CrMo}) For any \'etale groupoid $\G$ of dimension $d$, one has
%\[ H^{n}(\G; \oo) \cong H^{d-n}_{c}(\G; \mathbb{R})^{\vee} \ .\]
%\end{cor}

%{\it Proof: } Recall that $\oo$ is defined by $\Gamma(U; \oo)= H^{d}_{c}(U)^{\vee}$. We choose $\U$ a basis by contractible
%opens. Since $\oo$ is locally constant, it follows that $H^*(\G; \oo)\cong \check{H}^{*}_{\U}(\G; \A)$ 
%is computed by the $\check{C}$ech-type complex
%with coefficients $U\mapsto H_{c}^{d}(U)^{\vee}$. We now apply the second part of Theorem \ref{theoreml}
%to $\A= \mathbb{R}$. To compute $\check{H}^{*}_{c}(\G;\mathbb{R})$, we use the resolution of $\mathbb{R}$ by 
%differential forms. Since $H^{l}_{c}(U)= 0$ unless $l= d$,
%the Corollary follows from the fact that the dual of the homology of a chain complex coincides with 
%the cohomology of the dual cochain complex (this is true because we work over $\mathbb{R}$).  \ \ $\Boxe$\\ 

\begin{num}\label{CWetalgr}{\bf Chern-Weil for \'etale groupoids: }\rm Clearly all the constructions of Section  \ref{classes}
apply to any \'etale groupoid $\G$, provided we use the $\check{C}$ech-De Rham complexes mentioned in \ref{CDRetale}.
 Hence, for any principal $G$-bundle
$P$ endowed with a smooth action of $\G$, one has an associated Chern-Weil map
\[ S(\mathfrak{g}^{*})^G \rmap H^{*}(\G; \mathbb{R}) \]
whose image vanishes in degrees $> 2d$, where $d= dim(\G)$. The refined characteristic map,
\[ H^{*}(W_{d}(\mathfrak{g}, K)) \rmap H^{*}(\G; \mathbb{R}) \]
defines the exotic characteristic classes. Of particular interest is the (frame bundle of the) tangent space of $\nG{0}$, which is naturally
endowed with an action of $\G$, and which induces the exotic characteristic map of $\G$,
\begin{equation}\label{exetale} 
k_{\G}: H^{*}(WO_{d}) \rmap \check{H}_{\U}^{*}(\G)\cong H^{*}(B\G; \mathbb{R}) \ .
\end{equation}
\end{num}

\hspace*{.1in}When $\G= Hol_{T}(M, \F)$ this is the map discussed in section \ref{classes}. But this is not the only interesting
example. For instance, if one works with foliated bundles which are not necessarily transversal
(as e.g. in \cite{KaTo}), then one has to replace the holonomy groupoid $Hol_{T}(M, \F)$ by the monodromy groupoid $Mon_{T}(M, \F)$. 
The new versions of Theorem \ref{theorem2} and Corollary \ref{corex} for foliated bundles then yield characteristic classes
in $H^*(BMon_{T}(M, \F))$. These classes are refinements of the characteristic classes in $H^*(M)$, already constructed in \cite{KaTo}. \\
\hspace*{.3in}Another interesting example is when $\G$ is Haefliger's $\Gamma^q$.
The importance of this example lies into the fact that $\Gamma^q$ plays a classifying role for codimension $q$ foliations,
hence its cohomology consists on ``universal'' classes. We will elaborate this in \ref{univfor} of the next section.\\

%***************************************************************
%***************************************************************
%***************************************************************
%***************************************************************
\section{Explicit formulas}
\label{explicit}
%***************************************************************
%***************************************************************
%***************************************************************
%***************************************************************

In this section we illustrate our constructions in the case of normal bundles.
In particular we deduce Bott's formulas for cocycles associated to group actions \cite{Boform},
as well as Thurston's formula.

\begin{num}{\bf Explicit formulas for the normal bundle: }\rm We now apply the construction
of the exotic characteristic map of Section \ref{classes} 
to the normal bundle $\nu$. Corollary \ref{corex} applied to the (principal $GL_q$-bundle associated to) $\nu$
provides us with a characteristic map
\[ k_{\F}: H^*(WO_q) \rmap \check{H}_{\U}^{*}(M/\F) \ ,\]
which, when composed with the pull-back $\pi^{*}: \check{H}_{\U}^*(M/\F)\rmap H^*(M)$, gives the familiar exotic characteristic classes \cite{Bott} of $\F$.
Here $WO_{q}$ is the standard \cite{Bott} simplification of the truncated relative Weil complex $W_q(\mathfrak{gl}_{q}, O(q))$ that we now recall.
The idea is that the relative Weil complex $W(\mathfrak{gl}_{q}, O(q))$ (see \ref{clCW})
is quasi-isomorphic to a smaller subcomplex, namely
the dg algebra $S[c_1,\, .\, .\, .\, , c_q] \otimes E(h_1, h_3, \, .\, .\, .\, , h_{2[\frac{q+1}{2}]-1})$
generated by elements $c_i$ of degree $2i$ (namely the polynomials $c_i(A)= Tr(A^i)$), elements $h_{2i+1}$ of degree
$4i+1$ (any elements which transgress $c_{2i+1}$), with the boundary
\[  d(c_i)= 0, \ \ d(h_{2i+1})= c_{2i+1} .\]
Truncating by polynomials of degree $> q$, the resulting inclusion
\[ WO_q:= S_q[c_1,\, .\, .\, .\, , c_q] \otimes E(h_1, h_3, \, .\, .\, .\, , h_{2[\frac{q+1}{2}]-1})\rmap W_q(\mathfrak{gl}_{q}, O(q))\]
induces isomorphism in cohomology. With this simplification, the desired cohomology can be computed explicitly. 
Apart from the classical Chern elements $c_i$ (non-trivial only for $i< q/2$ even) there are new exotic classes.
Referring to \cite{Godb} for the complete description of $H^*(WO_q)$, we recall here that the simplest such class
is the Godbillon-Vey class $gv= [h_{1}c_{1}^{q}]\in H^{2q+1}(WO_q)$.
We denote by $gv_{\F}\in \check{H}_{\U}^{*}(M/\F)$ the resulting cohomology class $k_{\F}(gv)$. Its pull-back to $H^{*}(M)$ is the usual 
Godbillon-Vey class of $\F$. More generally, the Bott-Godbillon-Vey classes $gv^{\alpha}= [u_1c_{\alpha_1} \ldots c_{\alpha_t}]$ 
(and their images $gv^{\alpha}_{\F}$) are defined for any partition $\alpha= (\alpha_1, \ldots , \alpha_t)$ of $q$ (i.e. $q= \sum \alpha_i$).\\
\hspace*{.3in}For explicit formulas, let us choose a basis $\U$ so that $\tilde{\U}$ 
are also domains of trivialization charts for $\nu$ (as in \ref{rmks2.1}).  
Let $J_{h}: U\rmap GL_q$ 
denote the Jacobian of $h: U\rmap V$ (any
holonomy embedding). Then the $J_{h}$'s are the associated transition functions of
the transversal bundle $\nu$. Locally, we choose the trivial
connection $\nabla_U$ over $U$. The corresponding $\nabla(h)$ are then given by the
connection $1$-forms:
\[ \omega_{h}:= J_{h}^{-1} dJ_{h} \in \Omega^1(U; \mathfrak{gl}_q),\]
for $h: U\rmap V$. We see that the Chern character $Ch_{\nu}\in \check{C}^2(\U,\Omega^*)$ is given by:
\[ (h_1, \ldots , h_p)\mapsto  (-1)^p\int_{t_0+ t_1+ \ldots + t_p\leq 1} exp(\  (t_1\omega_{h_1}+ t_2\omega_{h_2h_1}+ \ldots + t_p\omega_{h_p 
\ldots h_2h_1})^2\ ) dt_0dt_1 . . . dt_p\ .\]
For instance, the first class $C_1= ch_1(\nu)\in \check{C}^*(\U,\Omega^*)$ has the components
\[ C_1^{\ps 1, 1\pd}(h)= Tr( J_{h}^{-1}dJ_{h}), \ C_1^{\ps 0, 2\pd}= C_1^{\ps 2, 0\pd}= 0\ .\]
As we know, this class is cohomologically trivial. This can be seen directly, since $U_1\in \check{C}^1(\U,\Omega^*)$,
\[ U_{1}^{\ps 0, 1\pd}= 0, \ U_{1}^{\ps 1, 0\pd}(h)= log(\mid det(J_{h})\mid) \]
transgresses $C_1$. Computing the resulting closed cocycle $U_1C_{1}^{q}$ we see that
\end{num}

\begin{cor}\label{godvey} The Godbillon-Vey class $gv_{\F}\in\check{H}^{2q+1}(M/\F)$ is represented in the $\check{C}ech$- De Rham complex by the
cocycle $gv_{\F}$ living in bi-degree $(q+1, q)$:
\begin{equation}\label{forgv} 
gv_{\F}( h_1, \ldots , h_{q+1})= log(\mid
det(J_{h_1})\mid) h_{1}^* Tr(\omega_{h_2})
h_{1}^*h_{2}^*Tr(\omega_{h_3}) \ldots
h_{1}^* . . . h_{q}^* Tr(\omega_{h_{q+1}}). 
\end{equation}
\end{cor}

Similarly, computing $U_1C_{\alpha_1}\ldots C_{\alpha_t}$ for a partition $\alpha= (\alpha_1, \ldots , \alpha_t)$ of $q$,
we obtain the following formula, which explains Bott's definition of the cocycles associated to group actions \cite{Boform}.

\begin{cor}The Bott-Godbillon-Vey class $gv^{\alpha}_{\F}\in\check{H}^{2q+1}(M/\F)$ is
 represented in the $\check{C}ech$- DeRham complex by the closed cocycle $gv^{\alpha}_{\F}$ living in bi-degree $(q+1, q)$:
\[ gv^{\alpha}_{\F}( h_1, \ldots , h_{q+1})= log(\mid
det(J_{h_1})\mid)\cdot h_{1}^*\{
Tr[\ \omega_{h_2}\cdot h_{2}^*(\omega_{h_3})\cdot \ldots  (h_{\ps\alpha_{1}-1\pd}
\ldots h_2)^*(\omega_{h_{\alpha_1}})]\ \} \cdot \]
\[ (h_{\alpha_1}\ldots h_2 h_1)^*\{ Tr[ \ 
\omega_{h_{\ps\alpha_{1}+1\pd}}\cdot 
h_{\ps\alpha_{1}+1\pd}^*(h_{\ps\alpha_{1}+2\pd})\cdot  \cdots
(h_{\ps\alpha_{1}+\alpha_{2}- 1\pd} \cdots h_{\ps\alpha_{1}+1\pd})^*
(\omega_{h_{\alpha_2}})]\ \}\cdot  \cdots \]
\end{cor}

\begin{num}\label{univfor}{\bf Universal formulas.}\rm\ As pointed out in the previous section, the constructions
that we described for foliations apply to any \'etale groupoid. Due to its classifying properties,
the case of the Haefliger groupoid $\Gamma^q$ is of particular interest. We wish to explain how the $\check{C}$ech-De Rham 
model for $\Gamma^q$ can be used to derive, in an explicit and straightforward way, the known formulas and properties
of universal characteristic classes for codimension $q$ foliations. We emphasize that all these properties are now part of the 
folklore on characteristic classes
for foliations, but they are usually derived by non-trivial abstract arguments at the level of classifying spaces.\\
\hspace*{.3in}First of all we make a slight simplification of the $\check{C}$ech-De Rham complex of $\Gamma^q$. Choosing the basis $\U$ of $\mathbb{R}^q$
by discs, since any such disc is diffeomorphic to $\mathbb{R}^q$, we see that
the category $Emb_{\U}(\Gamma^q)$ is equivalent to the category which has only one object, and all the embeddings  $\mathbb{R}^q\rmap \mathbb{R}^q$
as arrows. 
Accordingly, we define $\check{C}(\Gamma^q; \Omega)$ as in the previous sections, except that we take products only over strings
\[ \mathbb{R}^q \stackrel{\sigma_1}{\rmap} \ldots \stackrel{\sigma_k}{\rmap} \mathbb{R}^q \]
of embeddings $\mathbb{R}^q\rmap \mathbb{R}^q$.
The main theorem of this section implies
\end{num}

\begin{cor} The $\check{C}$ech-De Rham complex $\check{C}(\Gamma^q; \Omega)$ computes $H^*(B\Gamma^q; \mathbb{R})$.
\end{cor}

Now we can describe the main (cohomological) universal properties of $\Gamma^q$ in an 
explicit (and obvious) fashion. 
First of all, the universal property of $\Gamma^q$ can be seen easily in cohomology: given any codimension $q$ foliation,
choosing a basis $\tilde{\U}$ for $M$ and a transversal basis $\U$ as in \ref{rmks2.1}, there is an obvious map
$\check{H}(\Gamma^q) \rmap \check{H}_{\tilde{\U}}(M)$, to be seen as the map induced in cohomology by the classifying map
$M\rmap B\Gamma^q$ of $\F$ (well defined up to homotopy). This map is the composition of the pull-back (\ref{zero})
with another obvious map
\begin{equation}\label{clasmap}
\check{H}^*(\Gamma^q)  \rmap \check{H}^{*}_{\U}(M/\F)  
\end{equation}
(compare to \cite{Bott}). Now, all the characteristic maps for codimension $q$ foliations are just the composition of the (\ref{clasmap})'s
with a universal map
\begin{equation}\label{univchar}
k_{q}: H^*(WO_q) \rmap \check{H}^*(\Gamma^q)\ .
\end{equation}
Again, with the $\check{C}$ech-De Rham complexes at hand this is obvious, and $k_{q}$ is not at all abstract: it is just the characteristic map
(\ref{exetale}) applied to $\G= \Gamma^q$ and can be described in terms of the trivial connection on $\mathbb{R}^q$ (compare to \cite{Bott}). In particular, all the formulas of section \ref{classes}
come from similar universal formulas in  $\check{C}(\Gamma^q; \Omega)$.\\
\hspace*{.3in}At the price of more complicated formulas, we can further simplify the complex $\check{C}(\Gamma^q; \Omega)$.
Indeed, since the cohomology of $\Omega^*(\mathbb{R}^q)$ is $\mathbb{R}$ concentrated in degree zero (Poincar\'e lemma),
we see that $\check{H}(\Gamma^q)$ is also computed by the $\check{C}$ech (subcomplex) with constant coefficients
\[ \check{C}(\Gamma^q):\ \ \  0\rmap \mathbb{R} \rmap \prod_{\mathbb{R}^q\stackrel{\sigma_1}{\rmap}\mathbb{R}^{q}} \mathbb{R} \rmap 
\prod_{\mathbb{R}^q\stackrel{\sigma_1}{\rmap}\mathbb{R}^{q}\stackrel{\sigma_2}{\rmap}\mathbb{R}^{q}} \mathbb{R} \rmap \ldots \]
To pass from $\check{C}(\Gamma^q; \Omega)$ to $\check{C}(\Gamma^q)$ one has to repeatedly apply the Poincar\'e lemma. After a
lengthy but straightforward computation (for the details see Lemma 3.3.8 in \cite{Crath}) we obtain:

\begin{lem}\label{uuv} An $n$-cocycle in the $C$ech-De Rham complex:
\[ u= u_0+ u_1+ \ldots + u_n, \ \ u_k\in \check{C}^k(\Gamma^q,\Omega^{n-k}) \]
represents the same class in $\check{H}^{n}(\Gamma^q)$ as the $n$-cocycle $\tilde{u}$  in the $\check{C}$ech complex $\check{C}^*(\Gamma^q)$, given by:
\[  \tilde{u}(\sigma_1, \ldots , \sigma_n)= \sum_{s=0}^{n} (-1)^{ n\ps s-1\pd + \frac{s\ps s-1\pd}{2}}\int_{I_{\sigma_1, \ldots , \sigma_s}} u_{n-s}(\sigma_{s+1}, \ldots , \sigma_n) .\]
Here, $I_{\sigma_1, \ldots , \sigma_s}$ is the $s$-cube:
\[ I_{\sigma_1, \ldots , \sigma_s}(t_1, \ldots , t_s)= \sigma_s(\sigma_{s-1}( \ldots \sigma_3(\sigma_2(\sigma_1(0)t_1)t_2)  \ldots ) t_{s-1}) t_s .\]
\end{lem}

If we apply this to the Godbillon-Vey cocycle (i.e. to the formula (\ref{forgv}) in the $\check{C}$ech-De Rham complex $\check{C}(\Gamma^q; \Omega)$),
we obtain the well-known Thurston's formula:

\begin{cor} The universal Godbillon-Vey class $GV\in H^3(B\Gamma^1)\cong \check{H}^{3}(\Gamma^1)$ is represented in $\check{C}(\Gamma^1)$  by the cocycle:
\[ \tilde{gv}_1(\sigma_1, \sigma_2, \sigma_3)= \int_{0}^{\sigma_1(0)}
 log(\mid \sigma_{2}'(t)\mid)
 \frac{\sigma_{3}''(\sigma_2(t))}{\sigma_{3}'(\sigma_2(t))}
 \sigma_{2}'(t) dt .\]
\end{cor}

%***************************************************************
%***************************************************************
%***************************************************************
%***************************************************************
\section{Relations to basic cohomology}
\label{secbasic}
%***************************************************************
%***************************************************************
%***************************************************************
%***************************************************************

In the previous sections we have seen various models for the cohomology of the leaf space,
all canonically isomorphic. Let us put
\begin{equation}\label{cohls}
H^{*}(M/\F)= H^{*}(Hol_{T}(M, \F)),\ \ H^{*}_{c}(M/\F)= H^{*}_{c}(Hol_{T}(M, \F))\ .
\end{equation}
The reader may choose one of the many models: Haefliger's model (as indicated  by the above notations) i.e. \ref{bar} 
applied to the holonomy groupoid reduced to any complete transversal $T$, the $\check{C}$ech-De Rham model that we have described
in section \ref{CDRcomplex} (cf. Proposition \ref{lema1}), or the classifying-space model (cf. Theorem \ref{theorem1}).
We emphasize however that the last model only works for the cohomology without restriction on the supports!\\
\hspace*{.3in}Here and in the next section we explain why these cohomology theories are suitable  
theories for the leaf space. We first compare them to the more familiar {\it basic cohomology}
(see e.g. \cite{minimal, Ser}), which is a different cohomology theory for leaf spaces.\\

\begin{num}{\bf Basic cohomology.}\rm\  Choosing a basis $\U$ of opens of a complete transversal $T$
(or any transversal basis for $\F$),
one defines $\Omega^{k}_{bas}(T/\F)$ as the cohomology of $\check{C}^{*}(\U, \Omega^{k})$ in degree $*= 0$.
This complex consists on $k$-forms on $T$ which are invariant under holonomy, hence it does not depend on
the choice of $T$ (up to canonical isomorphisms, of course). The resulting cohomology is denoted $H^{*}_{bas}(M/\F)$.
There is an obvious map (induced by an inclusion of complexes)
\begin{equation}\label{jmap} 
j: H^{*}_{bas}(M/\F)\rmap H^{*}(M/\F)\ .
\end{equation}
Similarly  one defines the basic cohomology with compact supports $H^{*}_{c, bas}(M/\F)$ \cite{minimal}. 
The corresponding complex $\Omega^{k}_{c, bas}(T/\F)$ is the homology of $\check{C}^{*}_{c}(\U, \Omega^{k})$ in 
degree $*= 0$, i.e., as in \cite{minimal}, the quotient of $\oplus_{U\in \U} \Omega^{k}_{c}(U)$ by the span
of elements of type $\omega- h^{*}\omega$ ($h: U\rmap V$ is a holonomy embedding, and $\omega\in \Omega_{c}^{k}(V)$).
Again, there is an obvious map
\begin{equation}\label{jcmap} 
j_{c}: H^{*}_{c}(M/\F)\rmap H^{*}_{c, bas}(M/\F)\ .
\end{equation}
In general, the maps (\ref{jmap}) and (\ref{jcmap}) are not isomorphisms. 
The basic cohomologies are much smaller then $H^{*}(M/\F)$; for instance 
$H^{*}_{bas}(M/\F)= 0$ in degrees $*>q$, and they are finite dimensional if $\F$ is riemannian and
$M$ is compact.  The price to pay is the failure of most of the familiar properties  
from algebraic topology (e.g., as discussed below, Poincare duality and characteristic classes).
However we point out that (\ref{jmap}) and (\ref{jcmap}) are isomorphisms when the naive leaf 
space is an orbifold. This was explained in 4.9 of \cite{CrMo}, but the reader should think
about the similar statement for actions of finite groups on manifolds, and the fact that the
cohomology (over $\mathbb{R}$) of finite groups is trivial. In particular (see also \cite{Mol}), we have \\

\begin{prop}\label{cptl} If $(M, \F)$ is a riemannian foliation with compact leaves, then
{\rm (\ref{jmap})} and {\rm (\ref{jcmap})} are isomorphisms.
\end{prop}

Another fundamental property of our cohomologies (\ref{cohls}) is

\begin{prop}\label{pdutr}(Poincar\'e duality) For any codimension $q$ foliation $(M, \F)$, 
\begin{equation}\label{pppo} 
H^{*}(M/\F; \oo) \cong H^{q-*}_{c}(M/\F)^{\vee} \ .
\end{equation}
\end{prop}

\hspace*{.1in}This (and the more general Verdier duality) has
been proved in \cite{CrMo}. 
 Note however that, with the $\check{C}$ech model in hand, the theorem becomes obvious.
This new proof of Poincar\'e duality can be viewed as a rather straightforward extension of the classical proof
for manifolds \cite{BoTu} (and can also be interpreted as the obvious duality
between the homology and the cohomology of the discrete category $Emb_{\U}(\G)$, cf. \ref{embdcat}).
In contrast, the basic cohomologies $H^{*}_{bas}(M/\F)$ and $H^{*}_{c, bas}(M/\F)$
satisfy Poincar\'e duality only in the riemannian case \cite{Ser}. 
In this case these dualities are compatible via (\ref{jmap}) and (\ref{jcmap}), and they coincide if 
the leaves are compact (see Proposition \ref{cptl}).
\end{num}

\begin{num}{\bf Characteristic classes.}\rm\ As we have seen, one of the main features of $H^*(M/\F)$ is that it 
contains the characteristic classes of the bundles over the leaf space (i.e. transversal bundles), and the Bott vanishing theorem and the construction of the exotic classes hold at this level. Regarding the groups $H^{*}_{bas}(M/\F)$, again, 
they are too small to contain these characteristic classes. But, as before, this is not seen in the case of 
riemannian foliations. The reason is that,
if $\F$ is riemannian, then the transversal metric induces a {\it transversal connection}, i.e. a connection which is invariant under holonomy.
Using this type of connections in the construction the characteristic maps $k_{\nu}$ of the normal bundle $\nu$, we see
that $k_{\nu}: S(gl_{q}^{*})^{inv}\rmap H^*(M/\F)$ vanishes in degrees $> q$. This stronger vanishing result
(at the level of $H^*(M)$), together with the construction of the refined exotic characteristic map, appears in \cite{LaPa}. 
Moreover,
using the explicit constructions of section \ref{classes}, we see that $k_{\nu}$ (and its exotic
versions) factors through the basic cohomology groups. This obviously applies to general transversal bundles. In 
conclusion,

\begin{prop} If $P$ is a transversal principal $G$-bundles over $(M, \F)$ which admits a transversal connection then the characteristic map $k_{P}: S(\mathfrak{g}^{*})^{G}\rmap H^*(M/\F)$ of $P$ (cf. Theorem {\rm \ref{theorem2}}) vanishes in 
degrees $> q$. Moreover the map $k_{P}$ (and its exotic version, cf. Corollary {\rm \ref{corex}}) factors through the basic cohomology groups:\rm
\[  \xymatrix{
S(\mathfrak{g})^{G} \ar@{.}[rrr] \ar[rrrd]^-{k_{P}} 
& & & H^{*}_{bas}(M/\F) \ar[d]^-{j} \\ 
& & & H^{*}(M/\F) 
} \]\end{prop}
\end{num}

%\[  \xymatrix{
%S(\mathfrak{g})^{G} \ar[rrr] \ar[rrrd]^-{k_{P}} \ar[rrrdd]_-{k_{P}}
%& & & H^{*}_{bas}(M/\F) \ar[d]^-{j} \\ 
%& & & H^{*}(M/\F) \ar[d]^-{\pi^*}\\
%& & & H^*(M)
%} \]

\begin{num}\label{intf}{\bf Integration along the leaves.}\rm\ Haefliger's original motivation \cite{minimal} for introducing 
$H^{*}_{c, bas}(M/\F)$ is the existence of an integration over the leaves map $\int_{\F}^{\,'}: H^{*}_{c}(M)\rmap 
H^{*-p}_{c, bas}(M/\F)$ when the bundle of vectors tangent to the leaves is oriented. We want to point
out the existence of a refined integration,
\begin{equation}\label{integr2}
\int_{\F}: H^{*}_{c}(M)\rmap H^{*-p}_{c}(M/\F)\ ,
\end{equation}
which, composed with the canonical map (\ref{jcmap}), gives precisely Haefliger's integration. 
Using the $\check{C}$ech model this map becomes obvious: choosing $\U$, $\tilde{\U}$ as in \ref{rmks2.1},
the integration over the plaques (with the induced orientation),
$\int: \Omega^{*}_{c}(\tilde{U})\rmap  \Omega^{*-p}_{c}(U)$, induces
a map at the level of the $\check{C}$ech-De Rham complexes associated to $\U$ and $\tilde{\U}$. \\
\hspace*{.3in}An alternative abstract definition of $\int_{\F}$ follows e.g. from the spectral sequences
of \cite{CrMo} by standard methods of algebraic topology (``integration over the fiber'' as an edge map). 
The Hochschild-Serre spectral sequence (i.e. Theorem 4.4 of \cite{CrMo} applied to $\pi: M\rmap M/\F$) 
takes the form $H^{s}_{c}(M/\F; \el^{t})\Longrightarrow H^{s+t}_{c}(M)$,
where $\el^{t}$ is a transversal sheaf whose stalk above a leaf $L$ is $H^{t}_{c}(\tilde{L})$.
 This second description provides us with qualitative information. E.g., if the holonomy covers of the
leaves are $k$-connected, we find that $\int_{\F}$ is isomorphism in degrees
$n-k\leq * \leq n$. Using Poincar\'e duality, it follows that the pull-back map $H^{*}(M/\F)\rmap H^*(M)$ is isomorphism in degrees $0\leq *\leq k$. 
\end{num}

%***************************************************************
%***************************************************************
%***************************************************************
%***************************************************************
\section{Relations to foliated cohomology}
%***************************************************************
%***************************************************************
%***************************************************************
%***************************************************************

Another standard cohomology theory in foliation theory is the
{\it foliated cohomology} of foliations (see e.g. \cite{Alv, Hei, KT, MoSo}).
In contrast to the other cohomologies that we have seen
so far (transversal cohomologies), the foliated cohomology contains a great deal of longitudinal information.
In this section we describe its relation to our $\check{C}$ech-De Rham cohomology.

\begin{num}{\bf Foliated cohomologies.}\rm\ The foliated cohomology $H^*(\F)$ is defined in analogy with the De Rham cohomology
of $M$, which we recover if $\F$ has only one leaf. The defining complex is $\Omega^*(M, \F)= \Gamma(\Lambda^*\F)$, with the boundary defined by the usual Koszul-formula
\begin{eqnarray}\label{differential}
d(\omega)(X_1, \ldots , X_{p+1}) & = & \sum_{i<j}
(-1)^{i+j-1}\omega([X_i, X_j], X_1, \ldots , \hat{X_i}, \ldots ,
\hat{X_j}, \ldots X_{p+1})) \nonumber \\
 & + & \sum_{i=1}^{p+1}(-1)^{i}
L_{X_i}(\omega(X_1, \ldots, \hat{X_i}, \ldots , X_{p+1})) .
\end{eqnarray}
Here $L_{X}(f)= X(f)$. For later reference, we note the existence of an obvious (restriction to $\F$)
\begin{equation}\label{restr}
r: H^*(M)\rmap H^*(\F)
\end{equation}
There is also a version with compact supports, as well as versions $H^*(\F; E)$ with coefficients in any transversal (or foliated) vector bundle $E$:
one uses $E$-valued forms on $\F$, and one replaces the $L_{X_i}$ in the previous formula, by the derivatives $\nabla_{X_i}$
w.r.t. the Koszul connection of $E$ (see \ref{trbd}).
\end{num}

\begin{num}{\bf Remarks.}\rm\ In \cite{MoSo}, the cohomology $H^*(\F)$ is called
``tangential cohomology'', and is denoted $H^{*}_{\tau}(M)$. 
The groups $H^*(\F; \nu)$ with coefficients in the normal bundle
(see \ref{trbd}) first appeared in \cite{Hei} in the study of deformations of foliations, while those with 
coefficients in the exterior powers $\Lambda\nu$
show up e.g. in the spectral sequence relating the foliated cohomology with De Rham cohomology \cite{Alv, KT}. The groups $H^*(\F; E)$ with general coefficients can also be viewed as an instance of algebroid cohomology \cite{McK}.
Regarding the characteristic classes, since the Bott connection (see \ref{trbd}) is flat, it follows that the
characteristic classes of the normal bundle are annihilated by $r$. This new vanishing result at the level of foliated 
cohomology produces new (``secondary'') classes, $u_{4k-1}(\F)\in H^{4k-1}(\F)$. These appear in \cite{Gold} and have been described in great detail in \cite{Crave} in the more general context of algebroids. In particular, these new classes come 
from the cohomology groups $H^*(M/\F; \Omega^{0}_{bas})$ (via the map (\ref{Phi}) below). Still related to \cite{Crave},
let us mention that if $\F$ is the foliation induced by a regular Poisson structure on $M$, then one has an induced 
foliated bundle $K$ (the kernel of the anchor map), and $H^2(\F; K)$ contains obstructions to the integrability of the 
Poisson structure.
\end{num}

\hspace*{.1in}As explained in \cite{MoSo} in the case of trivial coefficients, and in \cite{Hei} in the case of the normal bundle as coefficients, the foliated cohomology can be expressed as the cohomology of certain sheaves on $M$.
For general coefficients $E$ we consider the sheaf $\Gamma_{\nabla}(E)$ described in \ref{trsh}.
A version of Poincar\'e's lemma with parameters shows that $H^k(\F; E)= 0$ in degrees $k>0$ if $\F$ is the standard 
$p$-dimensional foliation of $M= \mathbb{R}^p\times \mathbb{R}^q$. Since always $\Gamma_{\nabla}(M; E)= H^{0}(\F; E)$,
we deduce that $U\mapsto \Omega^*(\F|_{\, U}; E|_{\, U})$ is a flabby resolution of $\Gamma_{\nabla}(E)$, hence

\begin{prop}\label{folsheaf} For any foliated vector bundle $E$ over $(M, \F)$, $H^*(\F; E)$ is isomorphic to
$H^{*}(M; \Gamma_{\nabla}(E))$, the cohomology of $M$ with coefficients in the sheaf of $\nabla$-constant sections 
of $E$. In particular, $H^{*}(\F)\cong H^{*}(M; \Omega_{bas}^{0})$. The same holds for compact supports.
\end{prop}

\begin{num}{\bf Comparison.}\rm\ We now note the existence of a canonical map
\begin{equation}\label{Phi} 
\Phi: H^{*}(M/\F; \Omega^{0}_{bas}) \rmap H^*(\F)\ .
\end{equation}
Recall that $\Omega^{0}_{bas}$, as a sheaf on $M$, is the sheaf of smooth function which are constant on the leaves.
This map has various interpretations. First of all, it can be viewed as a version
with coefficients of the pull-back map (\ref{map2.1}) (cf. also Proposition \ref{folsheaf}). Accordingly,
the simplest description is in terms of the $\check{C}$ech-De Rham model. Choosing $\U$ and $\tilde{\U}$
as in \ref{rmks2.1}, the left hand side of (\ref{Phi}) is computed by the cochain complex $\check{C}^*(\U; C^{\infty}(U))$,
which is obviously a subcomplex of the $t= 0$ column of 
$\check{C}^s(\tilde{\U}; \Omega^t(\tilde{U},\F)$. Now (\ref{Phi}) is the map induced in cohomology.
Alternatively, at least when the holonomy groupoid is Hausdorff, $H^{*}(M/\F; \Omega^{0}_{bas})$ coincide with the 
differentiable cohomology \cite{difcoh}
of the holonomy groupoid of $\F$, and $\Phi$ is precisely the associated Van Est map described in \cite{WeXu}.
It then follows from one of the main results of \cite{Crave} (applied to the holonomy groupoid) 
that $\Phi$ is an isomorphism in degrees $\leq k$ provided the leaves (or their holonomy covers) are $k$-connected.
As in the previous section (see \ref{intf}), the same result follows e.g. from the spectral sequences  of \cite{CrMo}.
\end{num}

\begin{num}\label{intfle}{\bf Integration along the leaves.}\rm\ If $\F$ is oriented, then we have an integration map
\begin{equation}\label{integr3}
\int_{\F}: H^{*}_{c}(\F)\rmap H^{*-p}_{c}(M/\F; \Omega^{0}_{bas})\ .
\end{equation} 
This map is dual to the Van Est map (\ref{Phi}) and can be viewed as a version of
the integration map (\ref{integr2}) with coefficients in the normal bundle (accordingly,
there are similar maps for any transversal vector bundle $E$ over $M$, cf. also \ref{trsh} and Proposition \ref{folsheaf}).
Again, as in the previous section (see \ref{intf}), this map (\ref{integr3}) becomes obvious
if one uses the $\check{C}$ech-De Rham model. \\
\hspace*{.3in}We want to point out here that the integration over the fibers that we have described 
clarifies the construction of the Ruelle-Sullivan current of a measured foliation (cf. e.g. Section 3 of \cite{CoOp}, or \cite{MoSo} p 126),
and also gives new qualitative information about it. Fix a transversal basis $\U$ for $\F$. A
smooth transverse measure $\mu$ is just a measure on each $U\in \U$, which is invariant w.r.t. holonomy embeddings.
Hence the integration against $\mu$ is simply a linear map  
\[ \int_{\mu}: \Omega^{0}_{c, bas}(M/\F)= H^{0}_{c}(M/\F; \Omega^{0}_{bas})\rmap \mathbb{R}\ .\]
Combining this with the integrations along the leaves (\ref{integr2}), (\ref{integr3}),
we can arrange our maps into a diagram
(see also (\ref{restr}))
\[  \xymatrix{
H^{p}_{c}(M) \ar[r]^-{\int_{\F}} \ar[d]_-{r} & H_{c}^{0}(M/\F) \ar[d] & \\
H^{p}_{c}(\F)\ar[r]^-{\int_{\F}} & H_{c}^{0}(M/\F; \Omega^{0}_{bas}) \ar[r]^-{\int_{\mu}} & \mathbb{R}
} \]
The resulting map $\int_{C}: H^{p}_{c}(M) \rmap \mathbb{R}$ is precisely the integration
of \cite{CoOp} against the Ruelle-Sullivan current $C= C_{\mu}$ (and this defines $C$ as a degree $p$ element in the 
closed homology of $M$). As pointed out in \cite{MoSo}, $C$ actually comes from the closed homology of $\F$.
In terms of our diagram this simply means that $\int_{C}$ factors through $H^{p}_{c}(\F)$. 
\end{num}

\begin{num}{\bf Spectral sequences.}\rm\ Almost all of the maps that we have described in the last two sections
figure in certain the spectral sequences.  First of all, the filtration on $\Omega^*(M)$ induced by $\F$ 
(cf. e.g. \cite{Alv, KT}) induces a spectral sequence
\[ E^{s, t}_{1}= H^{s}(\F; \Lambda^t\nu) \Longrightarrow H^{s+t}(M)\ .\]
Similarly, the filtration of the $\check{C}$ech-De Rham double complex induces a spectral sequence
\[ \bar{E}^{s, t}_{1}= H^{s}(M/\F; \Omega^{t}_{bas}) \Longrightarrow H^{s+t}(M/\F)\ .\]
Note that $E^{0, t}_{2}= H^{t}_{bas}(M/\F)$. These two spectral sequences are related by the pull-back map (\ref{map2.1}), and by the Van Est map (\ref{Phi})
with coefficients, $\Phi: H^{s}(M/\F; \Omega^t)\rmap  H^{s}(\F; \Lambda^t\nu)$. With the same arguments as above,
these maps are isomorphisms in degrees $0\leq s\leq k$, if the holonomy covers of the leaves are $k$-connected. The version
with compact supports of this discussion involves (\ref{integr2}) and the integrations $\int_{\F}: H^{s}_{c}(\F; \Lambda^t\nu)\rmap H^{s-p}_{c}(M/\F; \Omega^{t}_{bas})$ (cf. \ref{intfle} above).
\end{num}

\end{document}